\documentclass[preprint,11pt]{article}

\usepackage{graphicx,epsf,amsbsy,amsmath,amssymb}
\usepackage[latin1]{inputenc}

\font\bbfnt=msbm10
\def\bbR{\mbox{\bbfnt R}}

\newcommand{\mb}[1]{\mbox{\bfseries \itshape #1}}

\newcounter{nuthm}
\newenvironment{thm}{\refstepcounter{nuthm}
    \begin{trivlist}\item[\hskip\labelsep{\bf Theorem~\thenuthm.}]\sl}
    {\rm\end{trivlist}}

\newcounter{nuprop}

\begin{document}

\markboth{J.-M. Ginoux, C. Letellier \& L. O. Chua}{Topological analysis of
the memristive circuit}

\title{Topological analysis of chaotic solution \\
of a three-element memristive circuit}

\author{Jean-Marc Ginoux,\\ Laboratoire {\sc Protee}, I.U.T. de Toulon,\\ Université du Sud, BP 20132
F-83957 La Garde Cedex, France,\\Christophe Letellier,\\ CORIA UMR 6614, Universit\'e de Rouen,\\ BP 12 F-76801 Saint-Etienne du Rouvray cedex, France,\\ Leon O. Chua, \\Electronics Research Laboratory\\ and Department
of Electrical Engineering and Computer Sciences,\\ University of California,
Berkeley, CA 94720}

\maketitle

\begin{abstract}
The simplest electronic circuit with a memristor was recently proposed.
Chaotic attractors solution to this memristive circuit are topologically
characterized and compared to R\"ossler-like attrators.
\end{abstract}

\date{{\bf Keywords}: Memristor; electronic circuit; chaos; topological analysis.}

\section{Introduction}

In classical electronics, there are three passive circuit elements: the
resistor, the capacitor and the inductor. In 1971, Leon Chua introduced a
fourth ``missing'' element which he named a memristor --- for memory resistor
--- by using symmetry arguments \cite{Chu71}. Chua also derived the properties
of this element. However, it is only in 2008 that Strukov and co-workers
\cite{Str08} found a memristance arising in a nanoscale system in which
solid-state electronic and ionic transport are coupled under an external bias
voltage. In 1976 Chua and Kang generalized the memristor to a broader class
of nonlinear dynamical systems they called memristive systems \cite{Chu76}. Few chaotic memristive electronic circuits were already proposed \cite{Ito08} and \cite{Mut10} but they were four-dimensional. It is only recently that a
three-dimensional system was proposed to describe a memristive circuit
\cite{Mut10b}. This 3D system is the simplest three-element electronic circuit
producing chaotic behaviors since it is only made of two linear passive
energy-storage elements, and an active memristive device \cite{Chu76}. Such a
simple electronic circuit is not an algebraically minimal system for which only
five terms in the three right members of the set of ordinary differential
equations are allowed \cite{Zha97}. The system here studied  has five linear terms
and two nonlinear terms. The fact that this simple memristive circuit is not
a minimal system is an advantage since very often, minimal systems have very
tiny parameter domains associated with chaotic regimes and small attraction
basin. As a consequence, the simple memristive circuit has a quite large
attraction basin and a quite large domain of its parameter space over which
the system is chaotic.

In this letter, we will perform a topological analysis of chaotic attractors
of this memristive circuit. The subsequent part of this letter is
organized as follows. Section \ref{elecir} is devoted to the governing
equations and their geometric interpretation in terms of flow curvature
manifold. In section \ref{topan}  the topological anlaysis in terms of
branched manifold --- or template --- is presented. Section \ref{conc} gives
some concluding remarks.

\section{The governing equations and their geometric interpretation}
\label{elecir}

\subsection{The set of differential equations}

\begin{figure}[ht]
  \begin{center}
    \includegraphics[height=4.0cm]{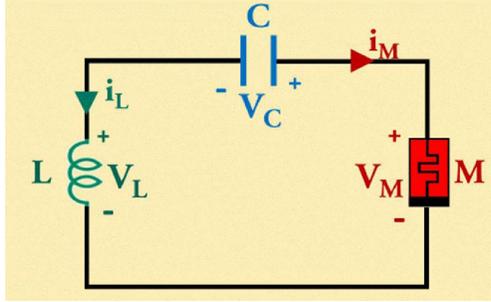} \\[-0.2cm]
    \caption{Block diagram of the simplest electronic circuit producing
chaotic behavior.}
    \label{memristor}
  \end{center}
\end{figure}

\begin{figure}[ht]
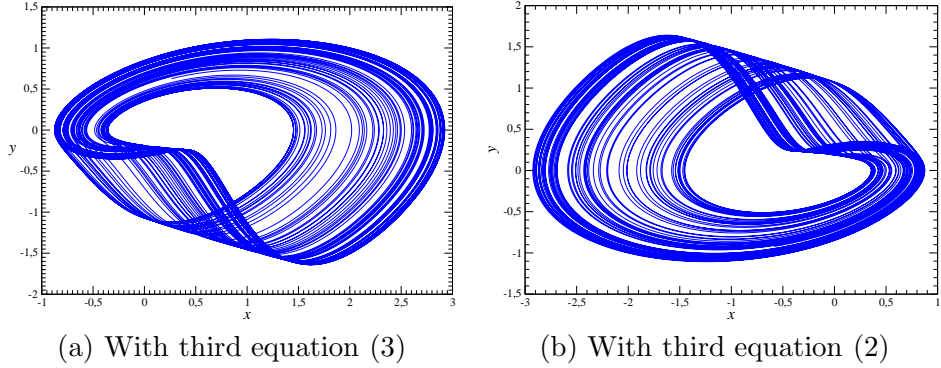

  \begin{center}
    \begin{tabular}{cc}
      \includegraphics[height=4.25cm]{mematt.eps} &
      \includegraphics[height=4.25cm]{memrimoins.eps} \\
      (a) With third equation (\ref{plus}) &
      (b) With third equation (\ref{moins}) \\[0.0cm]
    \end{tabular}
    \caption{Chaotic attractor solution to the minimal electronic circuit with
a memristor. Parameter value: $\alpha=0.98$.}
    \label{mematt}
  \end{center}
\end{figure}

The simplest electronic circuit producing chaotic attractors was proposed by
\cite{Mut10b}. The corresponding block diagram is shown in Fig.
\ref{memristor}. It is made of a linear passive inductor, a linear passive
capacitor and a nonlinear active memristor. This electronic circuit can be
described by the set of three differential equations
\begin{equation}
  \label{memeq}
  \left\{
    \begin{array}{l}
      \displaystyle
      \dot{x} = y \\[0.3cm]
      \displaystyle
      \dot{y} = - \frac{x}{3} + \frac{y}{2} - \frac{y z^2}{2} \\[0.3cm]
      \displaystyle
      \dot{z} = y - \alpha z - yz
    \end{array}
  \right. \, .
\end{equation}
It is written here in a slighly modified form since the third equation
is not \footnote{We picked a slightly different memristor characteristic
to demonstrate the simplest three-element chaotic circuit in Fig.
\ref{memristor} is robust with respect to the choice of the memristor
nonlinearity.}
\begin{equation}
  \label{moins}
  \dot{z} = - y (1-z) - \alpha z
\end{equation}
as in \cite{Mut10b} but rather
\begin{equation}
  \label{plus}
  \dot{z} = y (1-z) - \alpha z \, .
\end{equation}
The nonlinearity is thus inverted. As a consequence, the orientation of the
attractor in the $x$-$y$ plane is rotated by $\pi$ as easily checked in Fig.
\ref{mematt}. There is no difference observed in the topology of these two
attractors. This system has a single fixed point located at the origin of the
phase space.

\subsection{Flow curvature manifold}

Any trajectory solution to the dynamical system (\ref{memeq}) and denoted
$\mb{X} (t)$ describes the motion of a point M in the phase space
$\bbR^3(\mb{X})$ whose position depends on time $t$. This curve can also
be defined by its parametric representation
\begin{equation}
  \left\{
    \begin{array}{l}
      x = F_x (t) \\[0.1cm]
      y = F_y (t) \\[0.1cm]
      z = F_z (t)
    \end{array}
  \right.
\end{equation}
where $F_i$ are the right members of the governing equations (\ref{memeq}) that
are assumed to be continuous, $C^\infty $ functions. The {\it curvature}
$\kappa_1$ of any solution to system (\ref{memeq}) is defined as
\begin{equation}
  \label{curv1}
  \kappa_1 = \frac{ \mbox{Det} \left( \dot{\mb{X}}, \ddot{\mb{X}} \right)
  }{\left\| \dot{\mb{X}} \right\|^3} \, .
\end{equation}
Curvature measures the deviation of the curve from a straight line in the
neighborhood of any of its points.
In a similar way, the {\it torsion} $\kappa_2$ is defined as
\begin{equation}
  \label{curv2}
  \kappa _2 = \frac{ \mbox{Det}
  \left( \dot{\mb{X}}, \ddot{\mb{X}},\dddot{\mb{X}}  \right)
  }{\left\| \dot{\mb{X}} \wedge \ddot{\mb{X}} \right\|^2} \, .
\end{equation}
The torsion expresses the departure between the solution to system
(\ref{memeq}) and its projection into a place locally tangent to it.

It is known since Poincar\'e that the fixed points of a dynamical system
provides some information on the structure of the solution to that system
into the corresponding phase space. Fixed points belong to the zero-dimensional
invariat set of the dynamical system. It is also known that fixed points do
not provide all information needed to fully understand the structure of the
trajectories in the phase space. Good examples of such situation include
dynamical systems producing chaotic attractor structured around a single fixed
point. Few of them may be found in the collection of quadratic systems
proposed by \cite{Spr94}.

\newpage

To extract additional information, nullclines were sometimes used as in
\cite{Tho01}. But this is found only in the works made in the context of Fluid
Mechanics \cite{Rot98} where one-dimensional invariant sets were investigated
as connecting fixed points (when there are more than one) and structuring the
flow.
One-dimensional invariant set can be defined according to the following
theorem \cite{Gil10}.

\begin{thm}
The one-dimensional invariant set of a three-dimensional dynamical system is
defined as the location of the points where the velocity vector field is
co-linear to the acceleration vector field, that is, to the location of the
points where
\begin{equation}
  \label{eq9}
  \kappa_1 = 0 \quad
  \Leftrightarrow \quad \mbox{ Det} \left( {\dot{\mb{X}},\ddot{\mb{X}}}
  \right) = 0
\end{equation}

\end{thm}

More recently, one of us investigated two-dimensional invariant sets in the
context of the development of singularly perturbed systems or slow-fast
systems. In particular, a new geometric approach has been developed in order to
directly provide the slow invariant manifold of any $n$-dimensional dynamical
system, singularly perturbed or not \cite{Gin09}. Thus flow curvature manifold
can be defined as follows.

\begin{thm}
The flow curvature manifold of a three-dimensional dynamical system is defined
as the location of points where
\begin{equation}
  \label{eq10}
  \kappa_2 = 0 \quad
  \Leftrightarrow \quad
  \phi ( \mb{X} ) = \mbox{ Det}\left( {\dot{\mb{X}},\ddot{\mb{X}},
  \dddot{\mb{X}} } \right) = 0 \, .
\end{equation}
Solving this equation directly provides the analytical equation for the slow
invariant manifold associated with such a system.

\end{thm}

The flow curvature manifold is thus defined as \cite{Gin09b}
\begin{equation}
  \label{flocur}
  \phi (\mb{x}) = \dot{\mb{X}} \cdot ( \ddot{\mb{X}} \wedge \dddot{\mb{X}} ) = 0
  \, .
\end{equation}
Differentiating the acceleration vector $\ddot{\mb{X}}= \mb{J} \dot{\mb{X}}$
where $\mb{J}$ is the functional Jacobian matrix of the system, and inserting
this expression into (\ref{flocur}), we obtain
\begin{equation}
  \phi (\mb{X}) =
  \underbrace{\dot{\mb{X}} \cdot \left( {\cal J} \dot{\mb{X}} \wedge {\cal J}
  \ddot{\mb{X}} \right)}_{\phi_c} +
  \underbrace{\dot{\mb{X}} \cdot \left( \displaystyle \ddot{\mb{X}} \wedge
  \frac{{\rm d}{\cal J}}{{\rm d}t} \dot{\mb{X}} \right)}_{\phi_t}
\end{equation}
where $\phi_c$ is the time-independent component and $\phi_t$ is the
time-dependent component. Since $\phi_c$ does not contain the time derivative
of $\mb{J}$ it is associated with the linear component of the vector field
and $\phi_t$ with the nonlinear component. In the neighborhood of fixed points
$\mb{X}$, the time-independent component of the flow curvature manifold
corresponds to the osculating plane. As a consequence, the attractor takes the
shape of $\phi_c$ in this neighborhood because the osculating plane cannot be
crossed by a trajectory \cite{Gin09b}.

\begin{figure}[ht]
  \begin{center}
      \includegraphics[height=8.0cm]{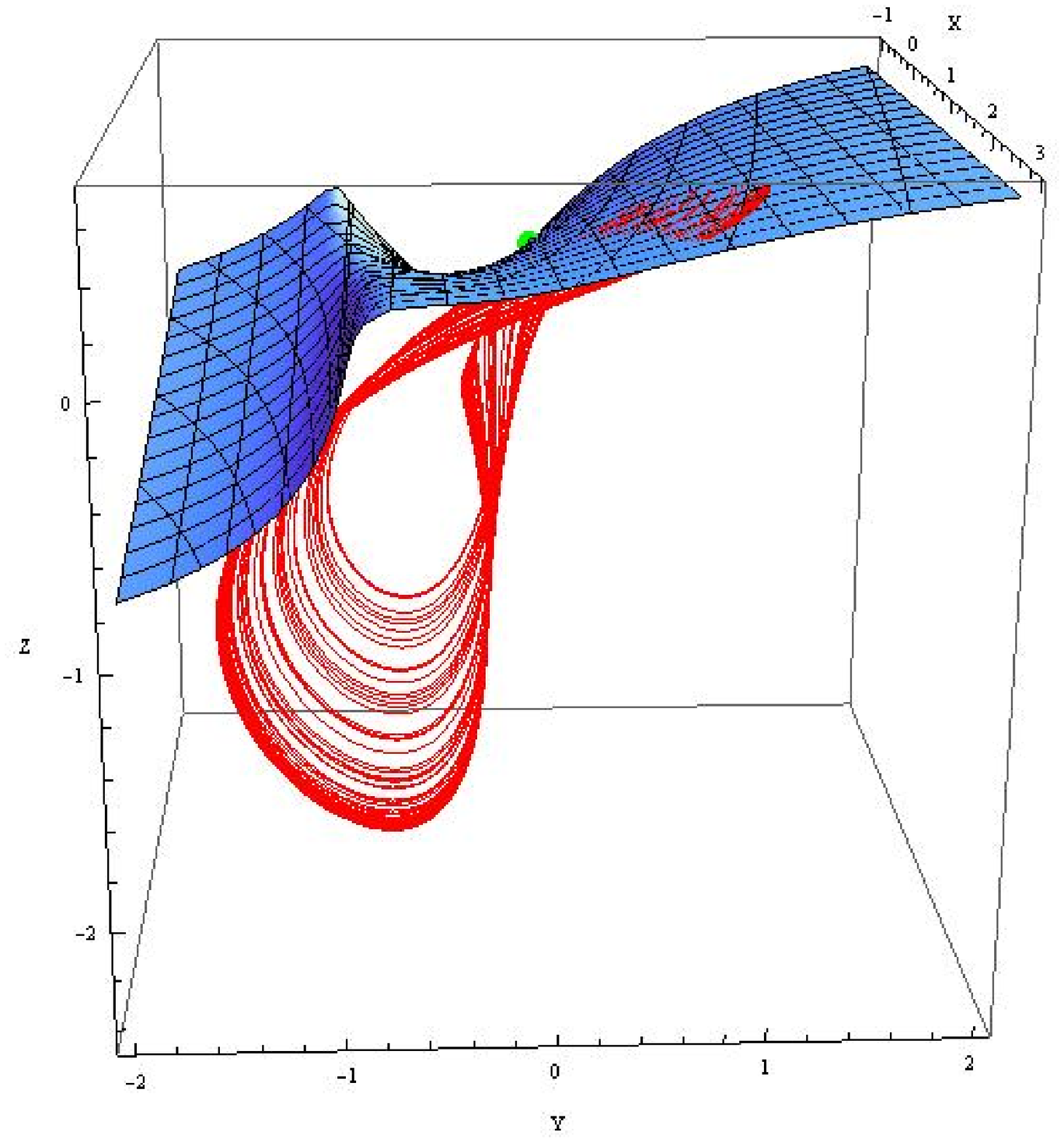}\\
      \includegraphics[height=8.0cm]{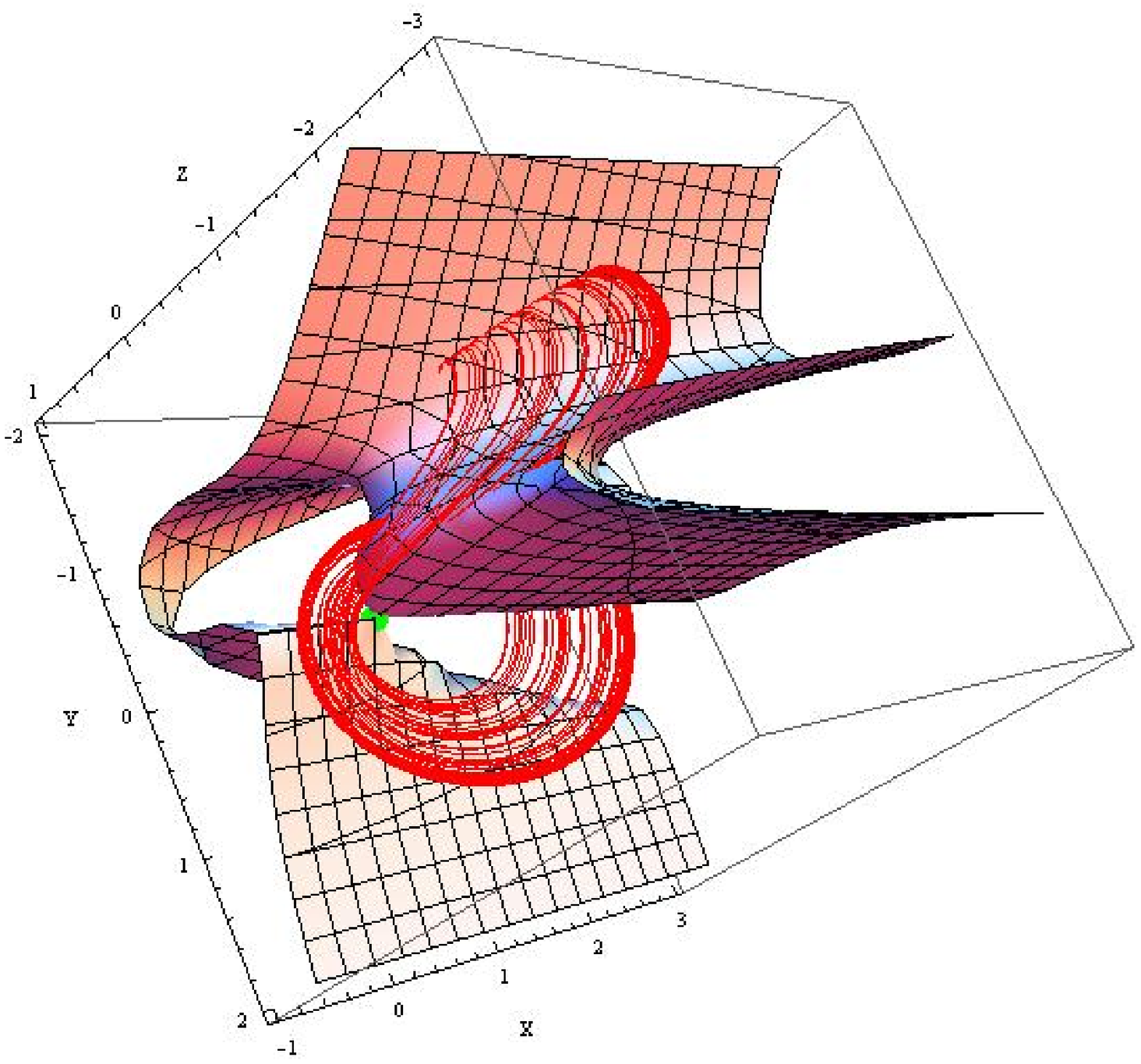}\\
      (a) Time-independent component $\phi_c$\\
      (b) Time-dependent component $\phi_t$ \\
    \caption{The two components of the flow curvature manifold $\phi$ for
the memristive circuit with parameter values: $\alpha=0.98$.}
    \label{memriflocur}
  \end{center}
\end{figure}

The two components of the flow curvature manifold of the memristive circuit
(\ref{memeq}) are shown in Fig.\ \ref{memriflocur}a. As expected, in the
neighborhood of the fixed point, the time-independent component of the flow
curvature manifold is tangent to the osculating plane. The component $\phi_t$
presents an hyperboloid (Fig.\ \ref{memriflocur}b). The trajectory wraps
around a significant part of this surface. Close to the fixed point, the
trajectory crosses component $\phi_t$ as observed in the R\"ossler system
\cite{Gin09b}.  Such an intersection between the trajectory and component
$\phi_t$ could be an explaination for the limitation to the development of the
dynamics.

\section{Topological analysis}
\label{topan}

\begin{figure}[ht]
  \begin{center}
    \includegraphics[height=6.5cm]{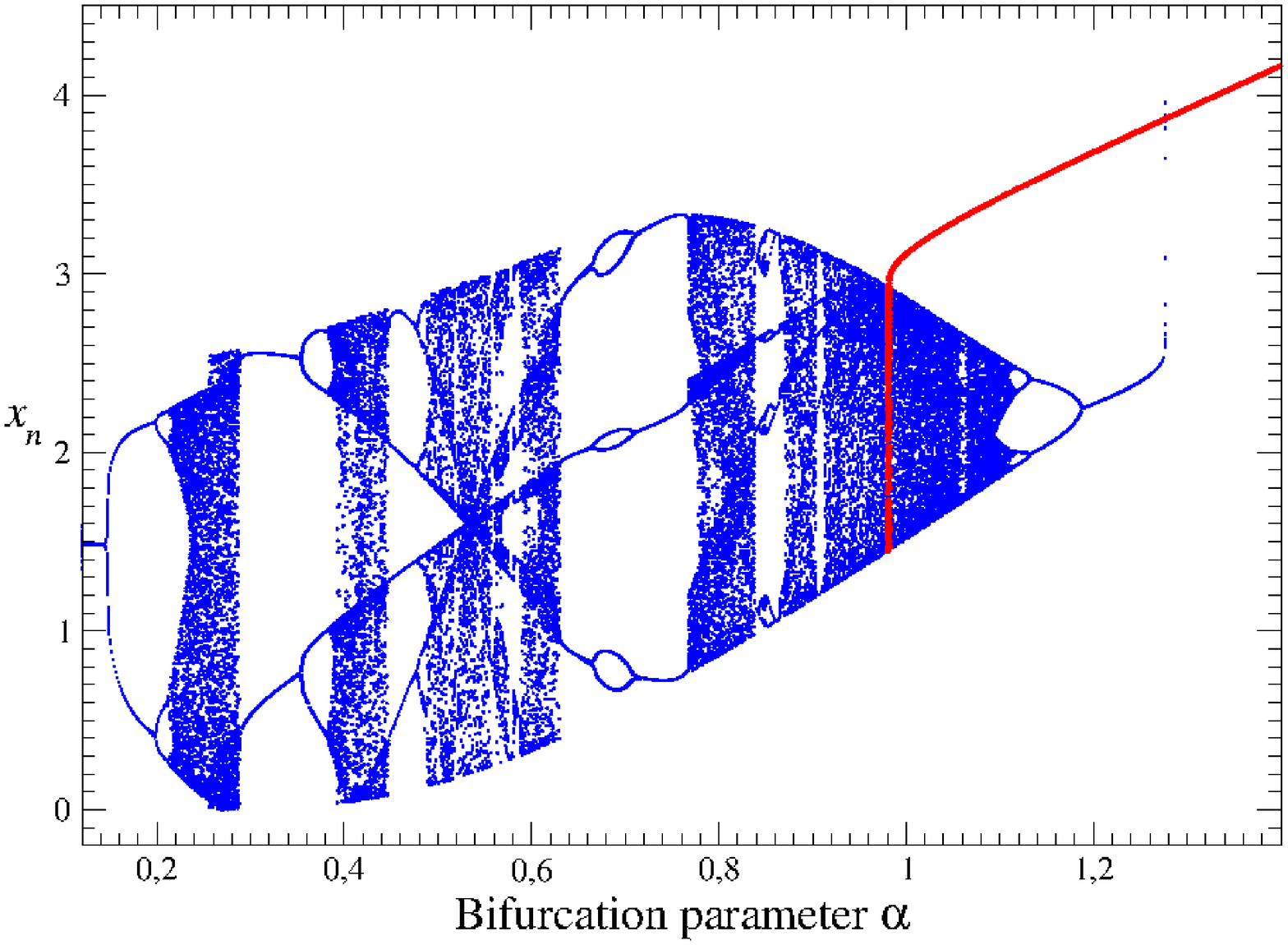} \\[-0.2cm]
    \caption{Bifurcation diagram versus parameter $\alpha$ of the minimal
electronic circuit with a memristor.}
    \label{memrybif}
  \end{center}
\end{figure}

For $\alpha=0.98$, the memristor electronic circuit (\ref{memeq}) produces
a chaotic attractor (Fig.\ \ref{mematt}). A Poincar\'e section can be safely
defined as
\begin{equation}
  \label{poisec}
  {\cal P} \equiv \left\{ x_n, z_n ~|~ y_n=0, \dot{y} < 0 \right\} \, .
\end{equation}
Using parameter $\alpha$ as a bifurcation parameter, a bifurcation diagram
is computed (Fig.\ \ref{memrybif}). When $\alpha$ is increased, the diagram
starts with a period-doubling cascade as a route to chaos. A crisis
is observed around $\alpha=0.25507$, a value at which the
attractor size suddenly increases. For slightly greater $\alpha$-values,
a second period-doubling cascade issued from a period-2 orbit is observed.
The development of the dynamics is observed at least up to $\alpha=0.53$.
Then for $\alpha$ slightly less than 0.6 a bubbling is observed.
The chaotic behavior is developed up to
another bubbling observed in the middle of the period-3 window
($\alpha \approx 0.69$). A last bubbling occurs in the middle of the second
period-3 window ($\alpha \approx 0.85$). For greater $\alpha$-value, the
dynamics is reduced by sub-critical bifurcations, destroying periodic orbits.
The diagram is ended by an inverse period-doubling cascade leading to a
period-1 limit cycle.

When $0.98 < \alpha < 1.276$, there are two attractors: one is chaotic and one
is a period-1 limit cycle. Depending on the initial conditions, the trajectory
thus settles down onto one of these two attractors. This bistability
disappears when the limit cycle collides with the chaotic attractor for
$\alpha$ is slightly greater than 0.98. We will thus investigate in details
the topology of chaotic attractors in the neighborhood of these two crises.

\begin{figure}[ht]
  \begin{center}
    \includegraphics[height=6.5cm]{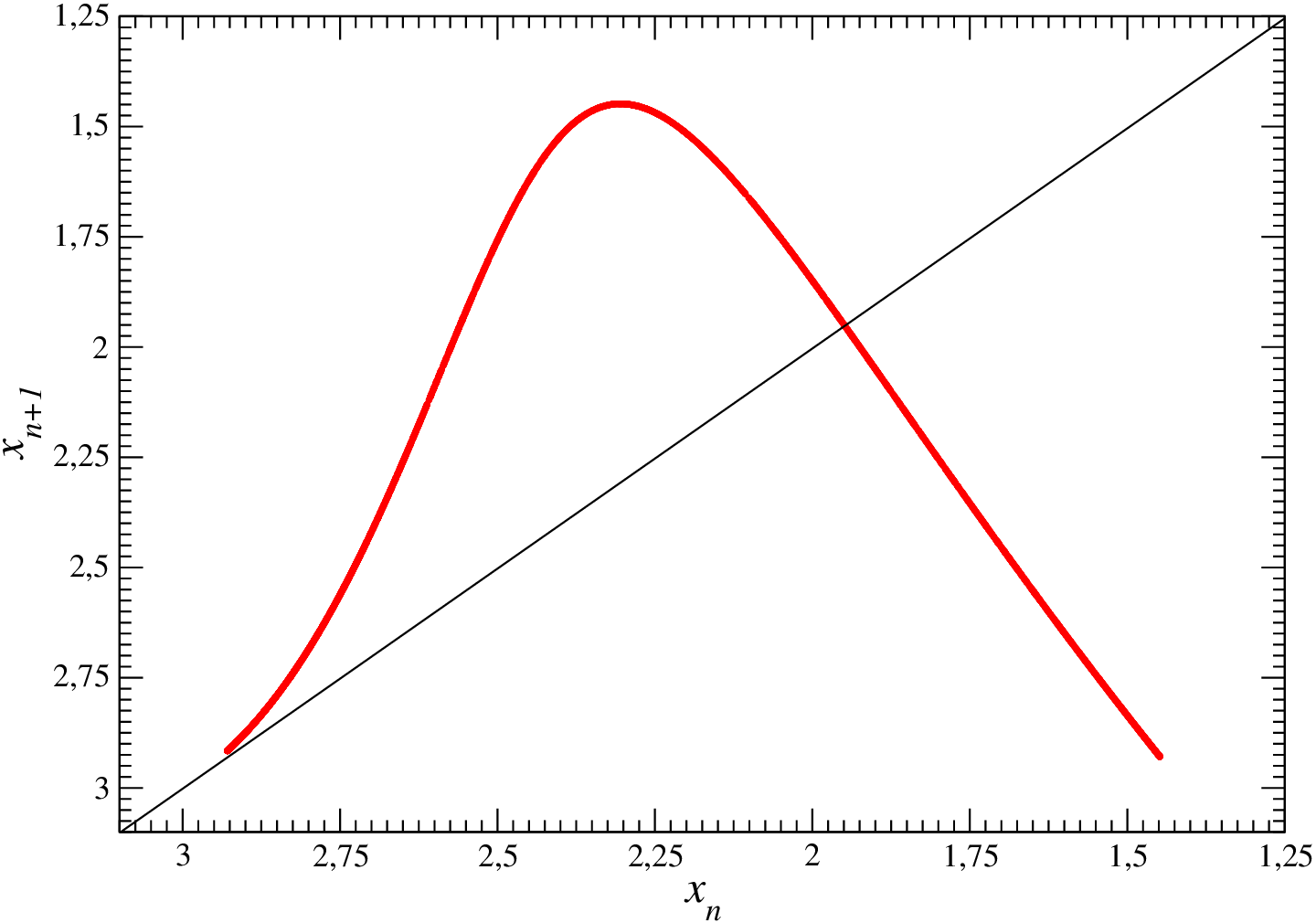} \\[-0.2cm]
    \caption{First-return map to a Poincar\'e section of the chaotic attractor
solution to the minimal electronic circuit with a memristor. Parameter value:
$\alpha=0.98$.}
    \label{memrymap}
  \end{center}
\end{figure}

Let us start with the chaotic attractor observed when the $\alpha$-value
is close to the boundary crisis, that is, when $\alpha=0.98$ (Fig.\
\ref{mematt}). A first-return map to the Poincar\'e section (\ref{poisec})
is computed (Fig.\ \ref{memrymap}). It is a unimodal map with a differential
maximum as expected after a period-doubling cascade. Such a maximum defines
a partition of the attractor that allows us to distinguish two topologically
distinct domains in the attractor \cite{Let95a}.
Increasing (decreasing) branches are  necessarily associated with domains ---
stripes --- with an odd (even) number of $\pi$-twists. In the present case,
the attractor is made of one stripe with one negative (anti-clockwise) and
$\pi$-twist and one stripe with two negative $\pi$-twists as detailed
below. We choose to label the increasing monotonic branch of the
first-return map by symbol ``2'' and the decreasing monotonic branch by
symbol ``1''. Each revolution over the attractor --- around the fixed point ---
is thus encoded according to
\begin{equation}
  \left|
     \begin{array}{lcl}
       1 & \mbox{ if } x_n > x_{\mbox{\footnotesize max}} \\[0.1cm]
       2 & \mbox{ if } x_n < x_{\mbox{\footnotesize max}}
     \end{array}
  \right.
\end{equation}
where $x_{\mbox{\footnotesize max}} = 2.305$. Using this encoding, a trajectory
in the phase space $\bbR^3 (x,y,z)$ is mapped into a symbolic sequence
\cite{Let95a}.

The boundary crisis occurs when the symbolic dynamics is complete, that is,
all symbolic sequences built with the two symbols ``1'' and ``2'' are realized
as unstable periodic orbits. This always arises when the end of one of the
monotonic branches touches the first bisecting line. In the present case (Fig.\
\ref{memrymap}) this is the increasing branch that touches the bisecting line.
The fact that a boundary crisis occurs exactly when the symbolic dynamics
is complete is a rather common feature since most of the quadratic minimal
systems investigated by \cite{Mal10} exhibit such property.

\begin{figure}[ht]
  \begin{center}
    \includegraphics[height=6.5cm]{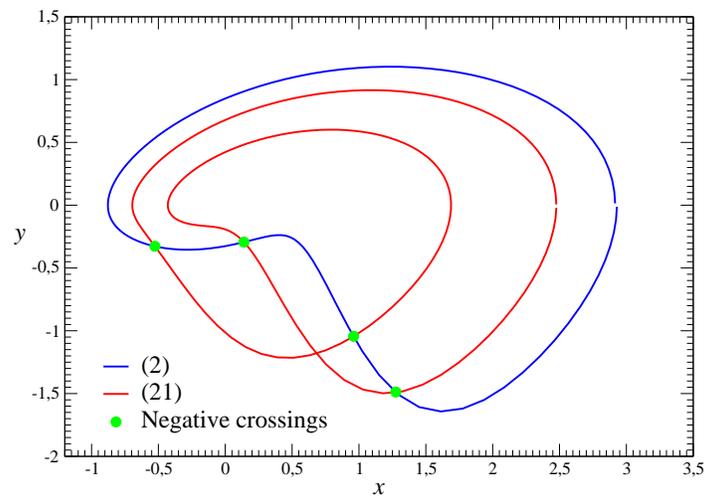} \\[-0.2cm]
    \caption{Pair of unstable periodic orbits extracted from the chaotic
attractor solution to the simplest electronic circuit with a memristor.
Parameter value: $\alpha=0.98$. The linking number Lk(2,21)=-2 according
to the four negative crossings identified in the $x$-$y$ plane projection.}
    \label{upo2-20}
  \end{center}
\end{figure}

The next step in the topological analysis is to extract the smallest unstable
periodic orbits and to compute some linking numbers. A linking number
counts the number of times a periodic orbit turns around another one. It can
be counted in a regular plane projection by identifying orientated crossings.
A crossing is counted only when one orbit ``crosses'' the other in the plane
projection (self-crossings are ignored). Crossings are then orientated
according to the third coordinates (see \cite{Let95a} for details). For
instance, orbits (2) and (21) exhibit four negative crossings (Fig.\
\ref{upo2-20}). The corresponding linking number is thus
\begin{equation}
  \mbox{Lk } (21,2) = \frac{1}{2} \left[ \displaystyle -4 \right] = -2 \, .
\end{equation}
This means that orbit (21) turns twice in the anti-clockwise direction around
orbit (2).

\begin{figure}[ht]
  \begin{center}
    \begin{tabular}{ccc}
      \includegraphics[height=7.5cm]{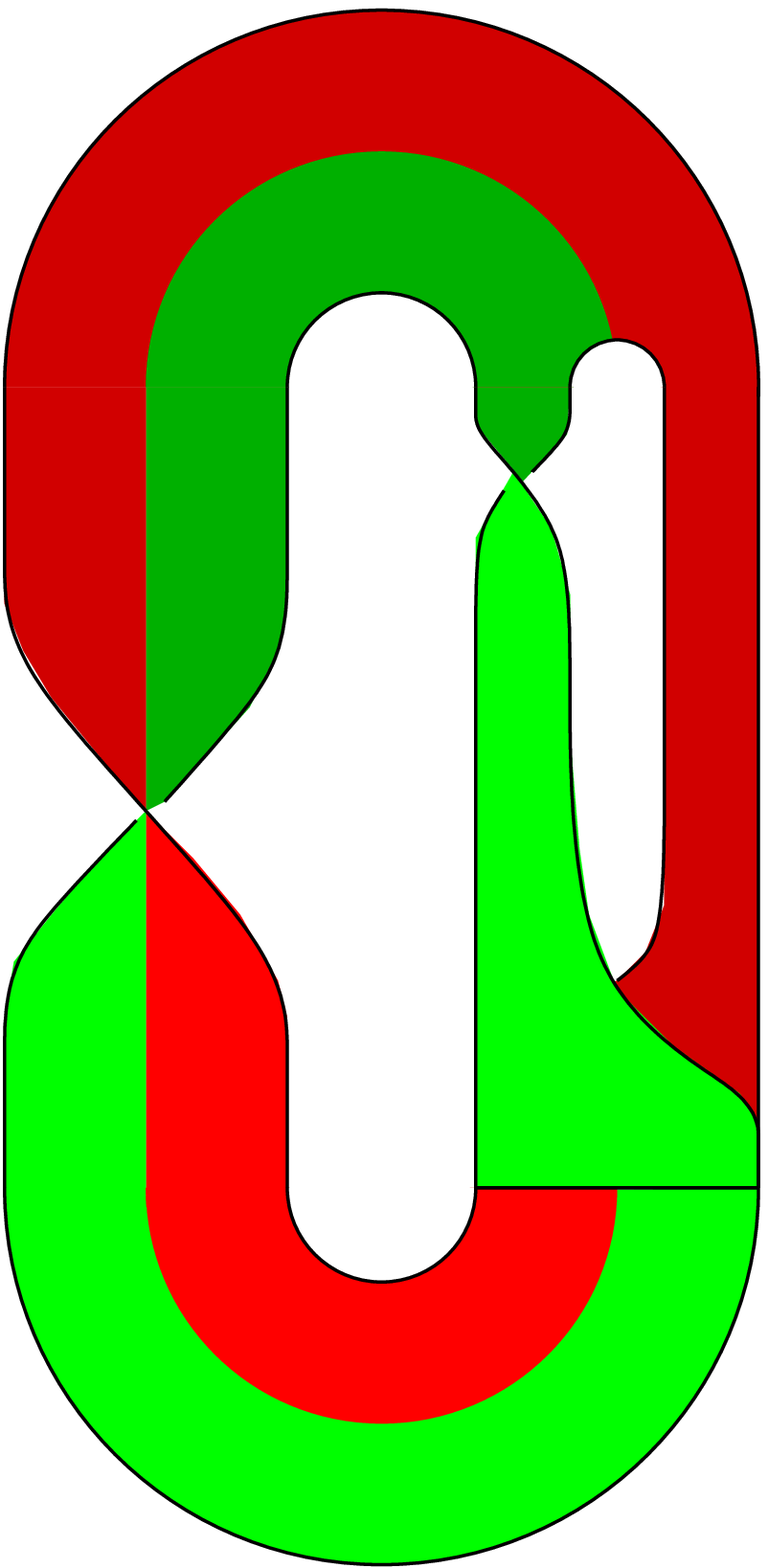} & ~ &
      \includegraphics[height=7.5cm]{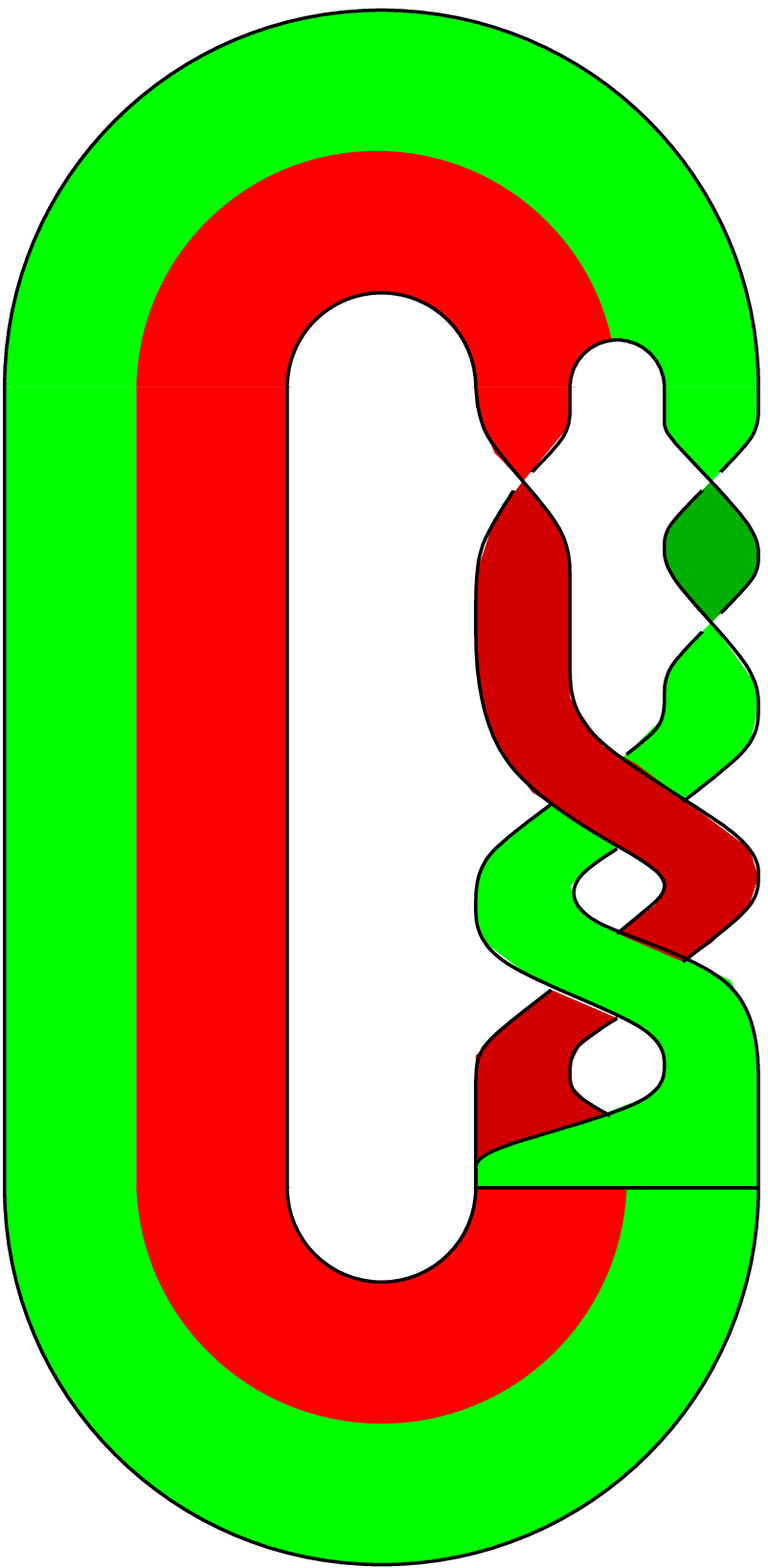} \\
      (a) Direct template & & (b) Standard template \\[-0.1cm]
    \end{tabular}
    \caption{Template for the chaotic attractor solution to the minimal
electronic circuit with a memristor. Parameter value: $\alpha=0.98$. The
global negative $\pi$-twist (a) is send to the right part (b) under an
isotopy.}
    \label{template}
  \end{center}
\end{figure}

We computed linking numbers for few pairs of periodic-orbits and checked
that all of them were correctly predicted by the template shown in Fig.\
\ref{template}a). The template exhibits a global negative $\pi$-twist (left
part of the template shown in Fig.\ \ref{template}a). As a consequence,
the attractor is a non trivial suspension of a unimodal map (a trivial
suspension would not present the global $\pi$-twist). The odd branches is
thus in the middle of the attractor and not at its periphery. Due to that, it
is said that the attractor is governed by an inverted Horseshoe map
\cite{Gil98}. The template proposed in Fig. \ref{template}a can be deformed
under an isotopy (a continuous deformation without any cutting) into a
``standard'' template according to \cite{Tuf92}. This standard representation
allows us to describe the template in terms of the following linking matrix:
\begin{equation}
  \label{linmat}
  M_{ij} =
  \left[
    \begin{array}{cc}
       -1 & -2 \\[0.1cm]
       -2 & -2
     \end{array}
   \right] \, .
\end{equation}
Diagonal terms $M_{ii}$ provides the numbers of (local) $\pi$-twist of each
stripe. As shown in the standard template (Fig.\ \ref{template}b) there is
therefore one negative $\pi$-twist for the odd --- red --- stripe (associated
with the decreasing branch of the first-return map) and two negative
$\pi$-twists for the even --- green --- strip (corresponding to the increasing
branch). Off-diagonal terms
$M_{ij}$ ($i \neq j$) means that stripe 1 crosses stripe 2 once in the
negative way as observed in the right part of the template shown in Fig.\
\ref{template}b.  With this convention, linking numbers can be algebraically
predicted from the linking matrix (\ref{linmat}) and the orbital sequences
\cite{Les94}.

\begin{figure}[ht]
  \begin{center}
    \includegraphics[height=6.5cm]{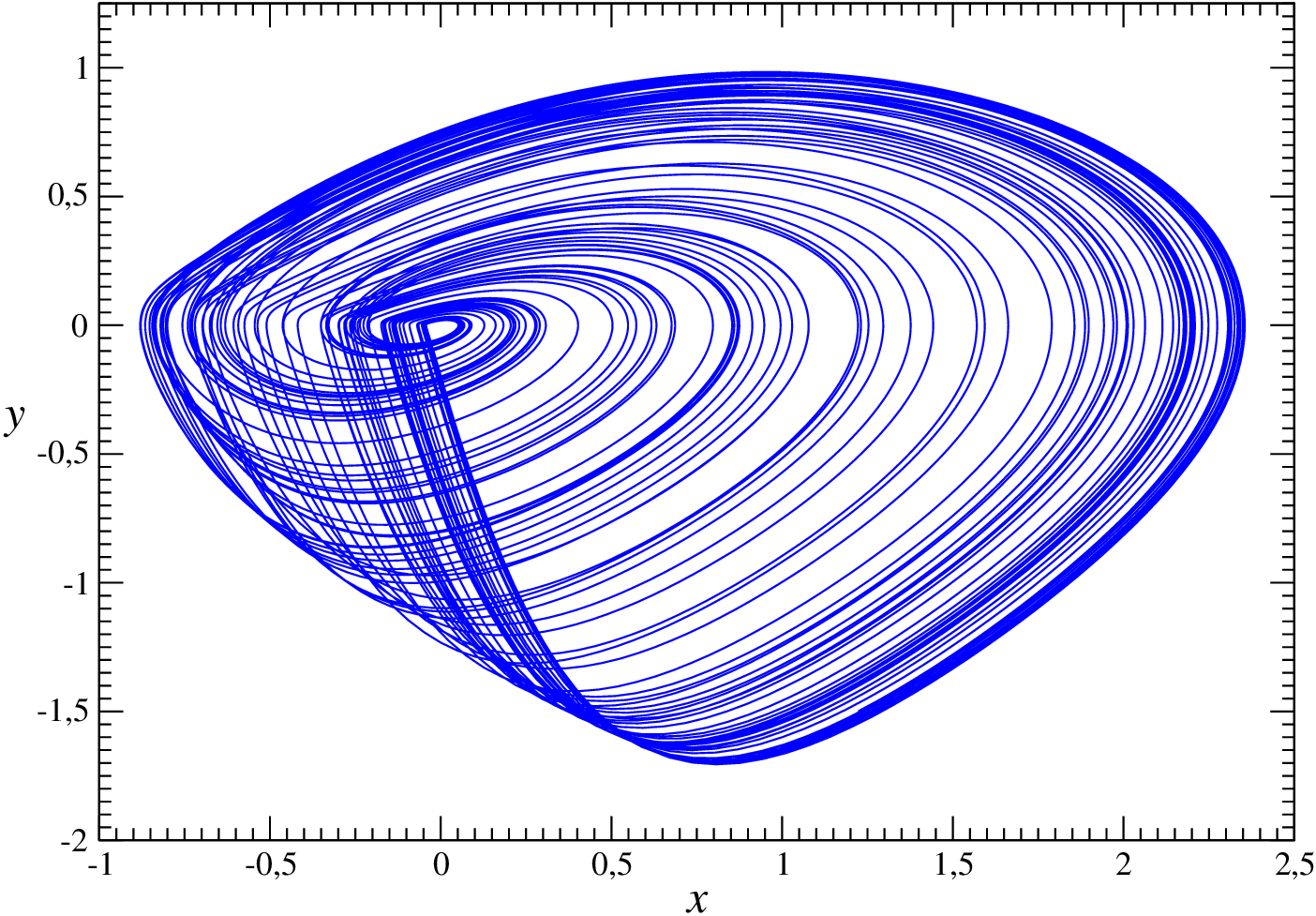} \\[-0.2cm]
    \caption{Chaotic attractor solution to the minimal electronic circuit with
a memristor. Parameter value: $\alpha=0.25507$.}
    \label{memrisde}
  \end{center}
\end{figure}

The second parameter value for which a topological analysis of the chaotic
attractor solution to system (\ref{memeq}) is chosen at the
crisis observed around $\alpha = 0.25507$. The chaotic attractor
(Fig.\ \ref{memrisde}) looks simpler than the attractor (Fig.\ \ref{mematt})
previously characterized. In particular, the global $\pi$ twist is no longer
seen.

\begin{figure}[ht]
  \begin{center}
    \includegraphics[height=6.5cm]{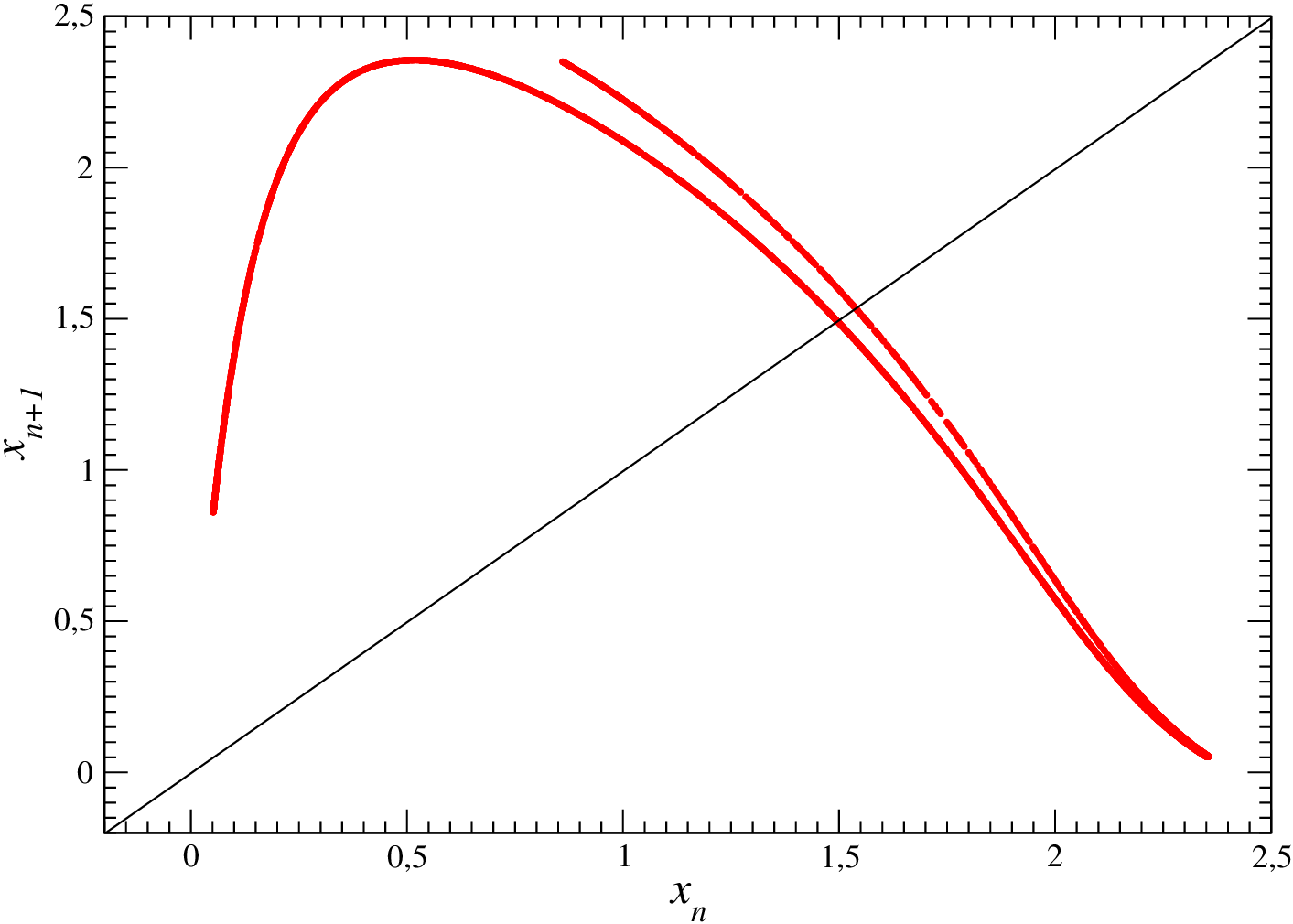} \\[-0.2cm]
    \caption{First-return map to a Poincar\'e section of the chaotic attractor
solution to the minimal electronic circuit with a memristor. Parameter value:
$\alpha=0.25507$.}
    \label{memrimade}
  \end{center}
\end{figure}

The first-return map to the Poincar\'e section (\ref{poisec}) confirms that
this is no longer an inverted Horseshoe since the increasing branch ---
associated with an even number of $\pi$-twists --- is now in the middle of
the attractor (lower $x$-value in Fig.\ \ref{memrimade}, thus closer to the
fixed point). This map has a ``layered'' decreasing branch. The end of the
upper branch is exactly at the same ordinate ($x_{n+1}$) as the differentiable
maximum. This special characteristic was sometimes encountered in few other
systems, and always associated with a crisis. This was for instance observed
in a minimal system proposed by \cite{Sch00}: in this case, a boundary crisis
that ejects the trajectory to infinity was observed \cite{Mal10}.

\begin{figure}[ht]
  \begin{center}
    \includegraphics[height=6.5cm]{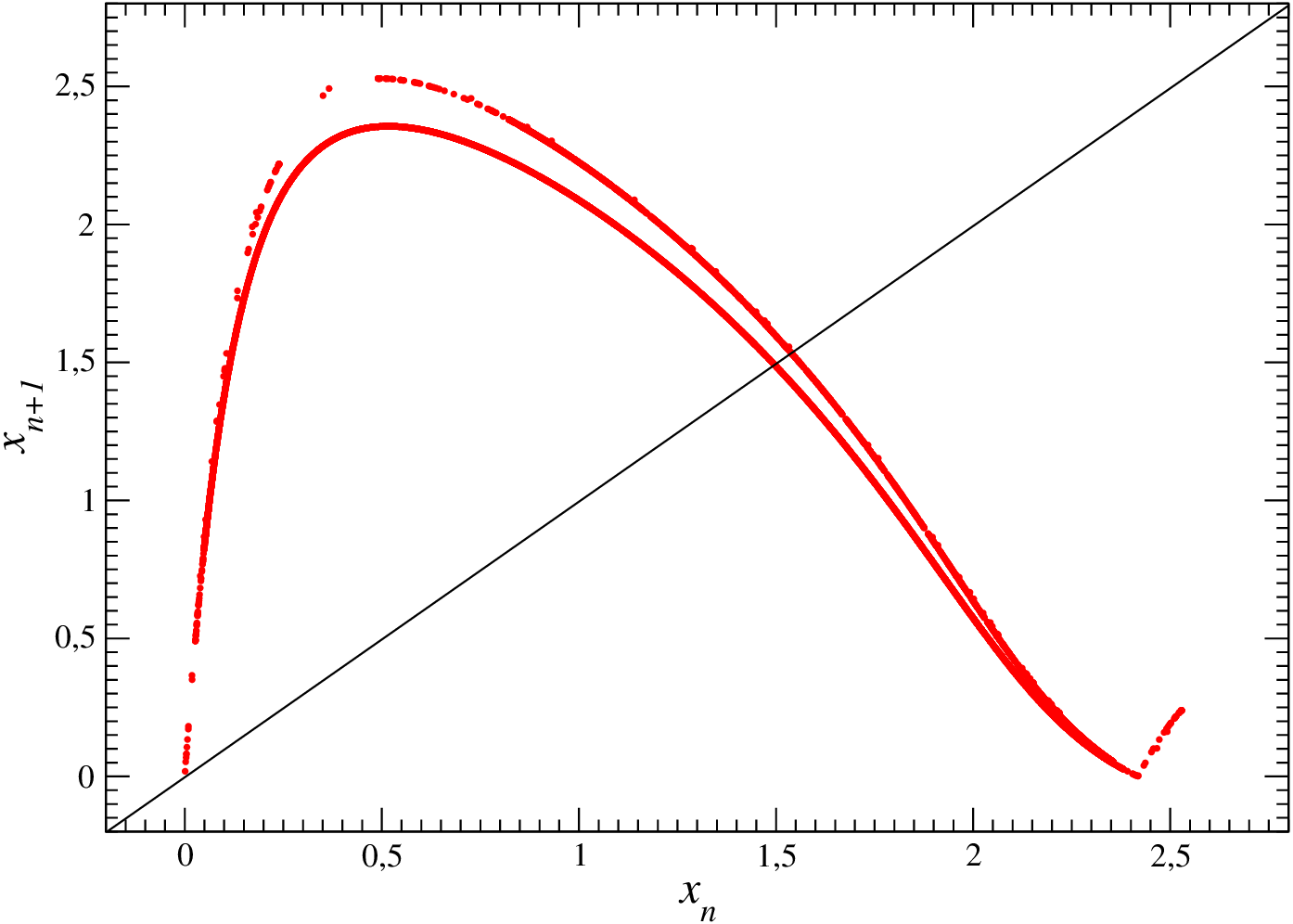} \\[-0.2cm]
    \caption{First-return map to a Poincar\'e section of the chaotic attractor
solution to the minimal electronic circuit with a memristor. Parameter value:
$\alpha=0.25515$.}
    \label{memap3}
  \end{center}
\end{figure}

The symbolic dynamics is obviously not complete (the increasing branch does
not touch the bisecting line) but as soon as the upper decreasing branch
has an ordinate greater than the maximum, there is an abrupt creation of
a large number of periodic orbits. Indeed, by increasing $\alpha$ to 0.25515,
the increasing branch is nearly completed and a third (increasing) branch
occurs (Fig.\ \ref{memap3}), thus revealing the explosive character of the
dynamics when the end of the upper decreasing branch has an ordinate that
becomes greater than the ordinate of the maximum (Fig.\ \ref{memap3}).

\begin{figure}[ht]
  \begin{center}
    \includegraphics[height=6.5cm]{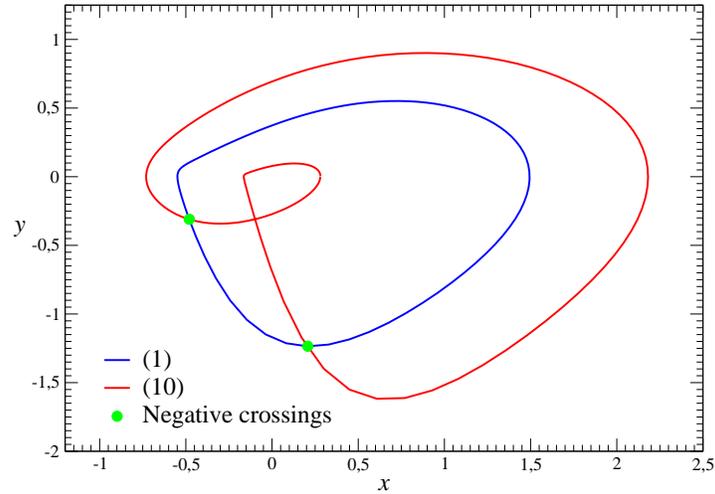} \\[-0.2cm]
    \caption{Pair of unstable periodic orbits extracted from the chaotic
attractor solution to the minimal electronic circuit with a memristor.
Parameter value: $\alpha=0.25507$. The linking number Lk(1,10)=-1 according
to the two negative crossings identified in the $x$-$y$ plane projection.}
    \label{upo1-10}
  \end{center}
\end{figure}

\begin{figure}[ht]
  \begin{center}
    \begin{tabular}{ccc}
      \includegraphics[height=6.5cm]{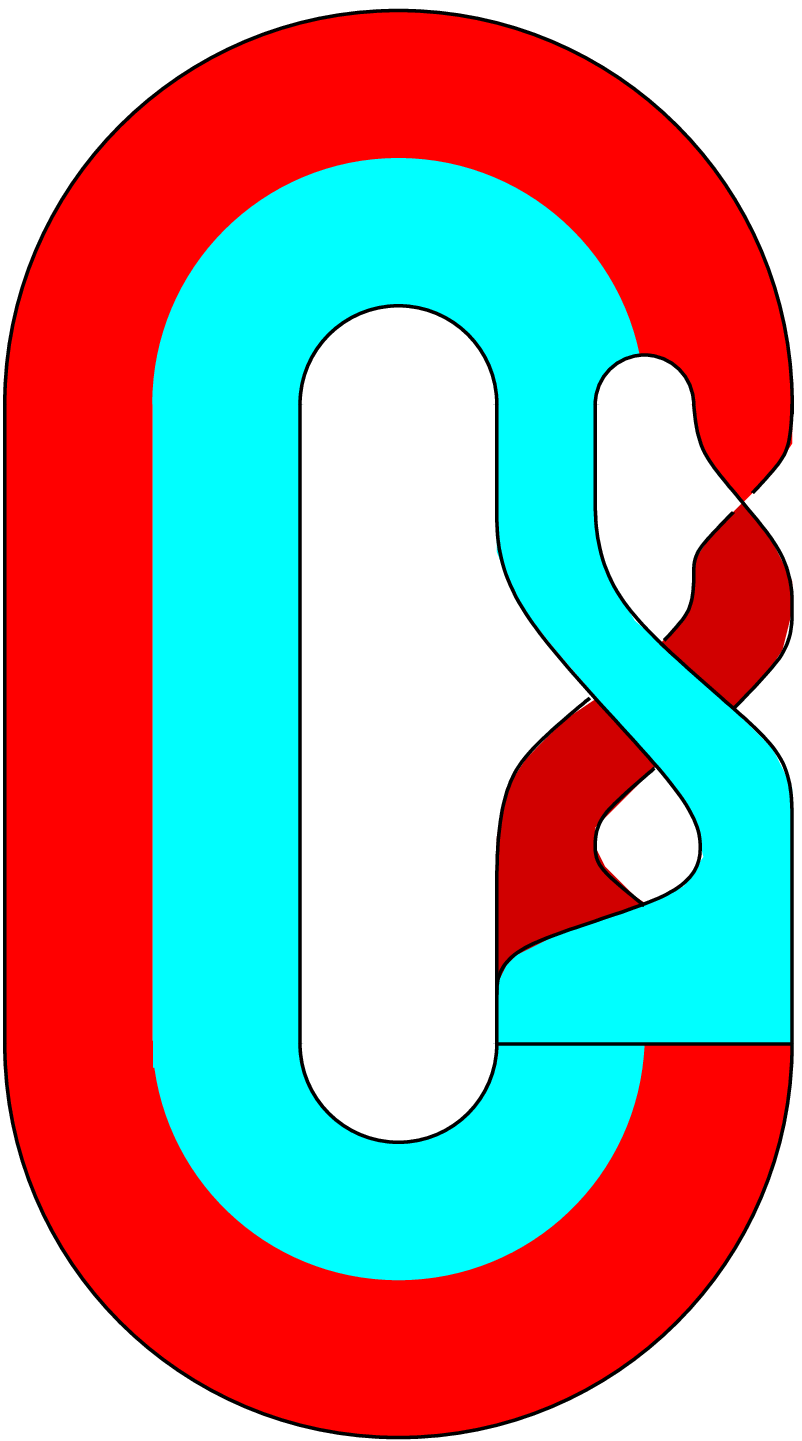} \\
    \end{tabular}
    \caption{Template for the chaotic attractor solution to the minimal
electronic circuit with a memristor. Parameter value: $\alpha=0.25507$.}
    \label{tempcris}
  \end{center}
\end{figure}

Periodic orbits were extracted for $\alpha=0.25507$. Linking numbers were
computed as shown in Fig. \ref{upo1-10}. All of the identified oriented
crossings found were negative. As expected after the first sight of the
chaotic attractor (Fig.\ \ref{memrisde}), the global negative $\pi$-twist
is no longer observed and all of the computed linking numbers were well
predicted by the template proposed in Fig.\ \ref{tempcris}. This is now
a trivial suspension of a unimodal map. The inner stripe --- close to the
fixed point --- has no local torsion and the outer stripe has a negative
$\pi$ twist. The linking matrix is therefore as
\begin{equation}
  M_{ij} =
  \left[
    \begin{array}{ccc}
      0 & -1 \\[0.1cm]
      -1 & -1
    \end{array}
  \right] \, .
\end{equation}
The cyan stripe associated with the increasing branch of the first-return map is
thus encoded with integer ``0'' and the decreasing red branch with integer
``1''.  Stripe 1 of this template is topologically equivalent to stripe 1 of
the template proposed for $\alpha=0.98$ (Fig.\ \ref{template}b). In fact,
this is the same stripe in both cases. When the bifurcation parameter $\alpha$
is varied from 0.25507 to 0.98, periodic points in stripe 0 are progressively
pruned while periodic points in stripe 2 are created.

\begin{figure}[ht]
  \begin{center}
    \includegraphics[height=6.5cm]{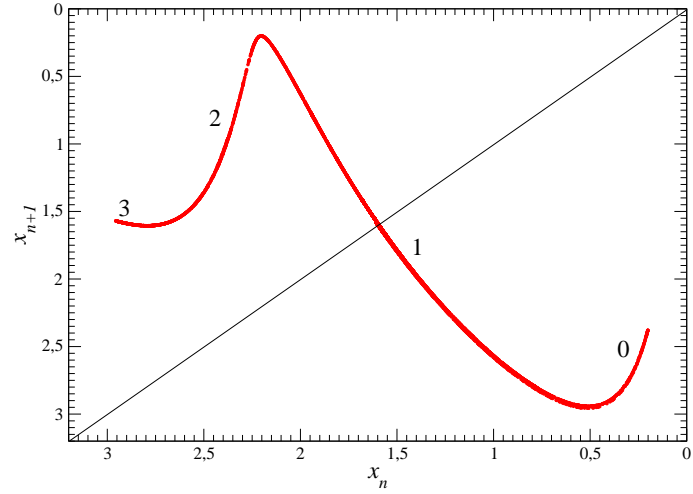} \\[-0.2cm]
    \caption{Three-modal first-return map to a Poincar\'e section of the
chaotic attractor solution to system (\ref{memeq}). Parameter values:
$\alpha=0.533$.}
    \label{memap4}
  \end{center}
\end{figure}

At intermediate value like $\alpha=0.533$, the three branches co-exist
(Fig. \ref{memap4}) in addition with a fourth one encoded by 4. It can be
shown in fact that the corresponding stripe present a local torsion equal
to four negative $\pi$-twists. The chaotic attractor can be therefore
split into four stripes whose topology can be synthetized by a template
described by the linking matrix
\begin{equation}
  M_{ij} =
  \left[
    \begin{array}{cccc}
      0 & -1 & -1 & -1 \\[0.1cm]
      -1 & -1 & -2 & -2 \\[0.1cm]
      -1 & -2 & -2 & -3 \\[0.1cm]
      -1 & -2 & -3 & -3
    \end{array}
  \right] \, .
\end{equation}
This matrix matches with those proposed for the R\"ossler attractor
of the funnel type \cite{Let95a}. From the topological point of view, the
dynamics produced by the memristive circuit is equivalent to those produced
by the R\"ossler system.

\begin{figure}[ht]
  \begin{center}
    \includegraphics[height=6.5cm]{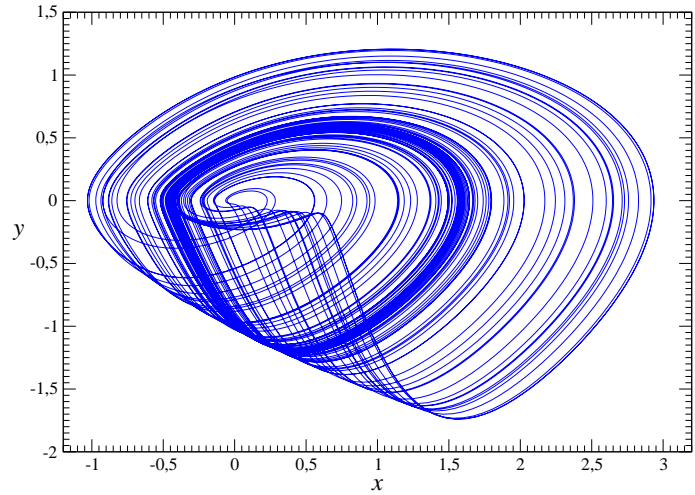} \\[-0.2cm]
    \caption{Chaotic attractor solution to system (\ref{memeq}). Parameter
values: $\alpha=0.533$.}
    \label{mematt4}
  \end{center}
\end{figure}

\begin{figure}[ht]
  \begin{center}
    \includegraphics[height=8.5cm]{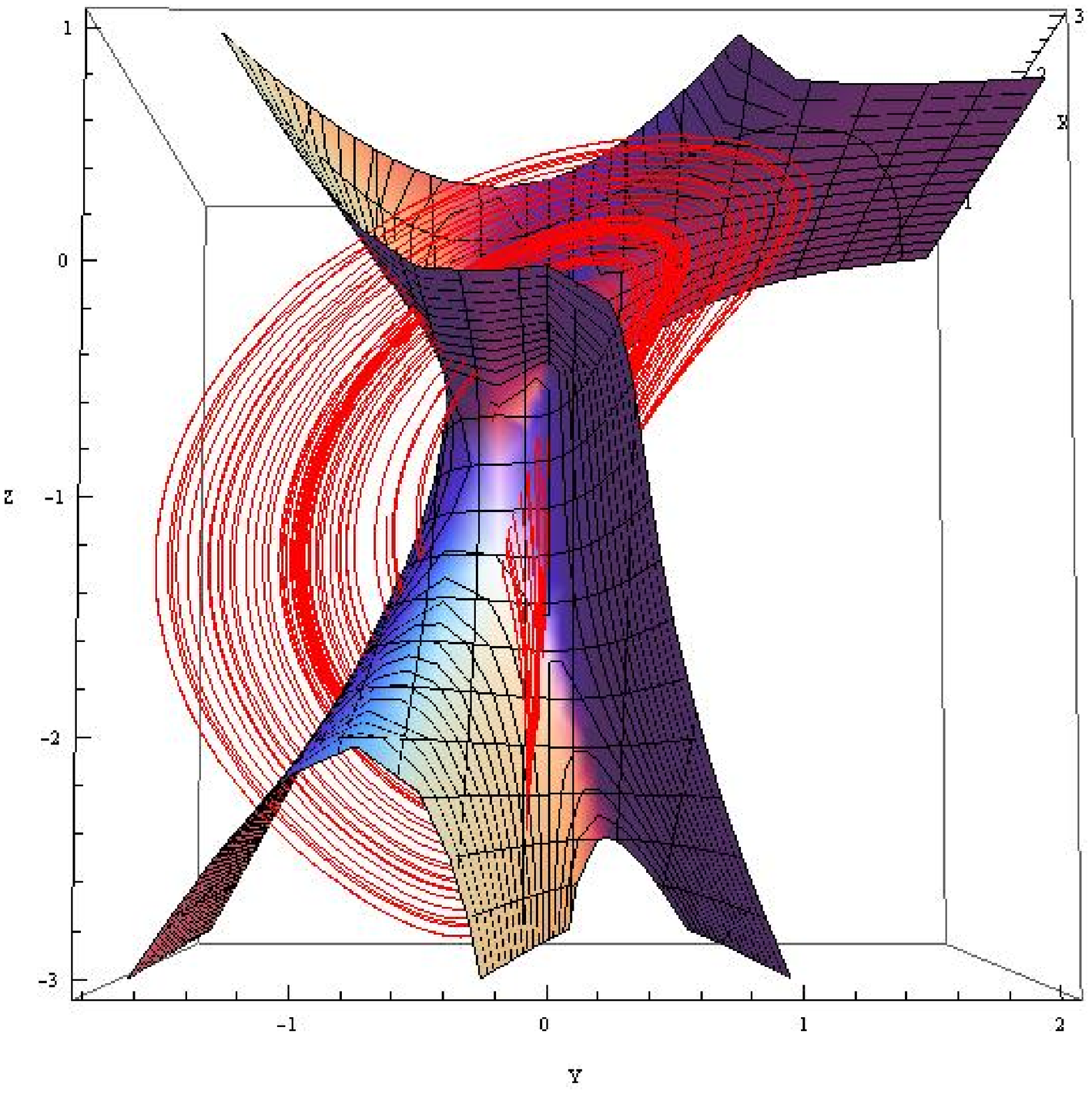} \\[-0.2cm]
    \caption{Time dependent component $\phi_t$ of the flow curvature manifold
for the memristive circuit. Parameter value: $\alpha=0.533$.}
    \label{timeflocur4}
  \end{center}
\end{figure}

The limitation of the dynamics could be explained by the time dependent
component of the flow curvature manifold $\phi_t$ as we did for few
R\"ossler-like systems in \cite{Gin09b}. From component $\phi_t$ as shown in
Fig.  \ref{timeflocur4}, it is obvious that the trajectory crosses component
$\phi_t$ in the neighborhood of the fixed point. Such a feature was always
remarked in different systems where some significant pruning was also noted.
This means that, for instance, when four branches are observed in the
first-return map (Fig.\ \ref{memap4}), all possible symbolic sequences are not
realized as unstable periodic orbits. Since crossings of the trajectory
occurs in the neighborhood of the fixed point and since this is the increasing
branch encoded by 0 that is the closes to the fixed point, symbolic sequences
with many consecutive ``0'' are forbidden. Consequently, the symoblic sequences
containing these sub-sequences are no longer realized as periodic orbits. This
is clearly seen in the four branches map (Fig.\ \ref{memap4}) where the
increasing branch ``0'' is almost removed.

\section{Conclusion}
\label{conc}

An electronic circuit with a memristor provides one of the simplest physical
realization of a dynamical system producing chaos. The governing equations
are not minimal from the algebraic point of view but this ensures in fact
the robustness of the chaotic behaviour so produced. Indeed, it is very often
observed that minimal systems have chaotic solutions for very tiny domain of
their parameter space and the attraction basin is quite limited. The simplest
memristive circuit is therefore a very good candidate for belonging to the
class of benchmark physical realizations for producing chaotic behaviors.

In this paper we investigate the topological structure of phase portraits
solution to this simple memristive circuit. It was shown that the chaotic
attractors produced by this circuit were of the R\"ossler-like type, that is,
they are topologically equivalent to the attractors produced by the well-known
R\"ossler system which results from a simplification of an oscillating
abstract chemical reaction. The dynamical regimes solution to the simplest
electronic circuit cannot be as developed in the R\"ossler system, due to a
strong limitation of the development of the dynamics induced by the time
dependent component of the flow curvature manifold. This leads to an inverse
period-doubling cascade that constrains
the behavior to be a limit-cycle for large value of the bifurcation
parameter.

\section*{Acknowledgements}
L. Chua acknowledged support from the U.S. Air Force. AFORS grant no. ??.\\
Authors would like to thank Dr. Ubiratan Freitas who made the cover figure.

\end{document}